# Probability models characterized by generalized reversed lack of memory property

**Asha Gopalakrishnan**[*]
**and Rejeesh C. John**

*Department of Statistics*
*Cochin University of Science and Technology*
*Cochin-22, Kerala, India.*
*e-mail:* asha@cusat.ac.in; rejeeshcjohn@yahoo.co.in

**Abstract:** A binary operator $*$ over real numbers is said to be associative if $(x*y)*z = x*(y*z)$ and is said to be reducible if $x*y = x*z$ or $y*w = z*w$ if and only if $z = y$. The operation is said to have an identity element $\tilde{e}$ if $x*\tilde{e} = x$. In this paper a characterization of a subclass of the reversed generalized Pareto distribution (Castillo and Hadi (1995)) in terms of the reversed lack of memory property (Asha and Rejeesh (2007)) is generalized using this binary operation and probability distributions are characterized using the same. This idea is further generalized to the bivariate case.

**AMS 2000 subject classifications:** Binary associative operation, Reversed lack of memory property, Bivariate distributions, Characterization..

## Contents

## 1. Introduction

A binary operation $*$ over real numbers is said to be associative if

$$(x*y)*z = x*(y*z) \tag{1.1}$$

for all real numbers $x, y, z$. The binary operation $*$ is said to be reducible if $x*y = x*z$ if and only if $y = z$ and if $y*w = z*w$ if and only if $y = z$. The general reducible continuous solution of the functional equation (1.1) is (Aczel, 1966, p. 254)

$$x*y = g^{-1}(g(x) + g(y)) \tag{1.2}$$

where $g(\cdot)$ is a continuous and strictly monotone function provided $x, y, x*y$ belong to a fixed interval $A$ in the real time. The function $g(\cdot)$ in (1.2) is determined up to a multiplicative constant; for all $x, y$ in a fixed interval $A$, implies $g_2 = \alpha g_1$ for all $x$ in that interval, for $\alpha \neq 0$. We assume hereafter that the binary operation is reducible and associative with the function $g(\cdot)$ continuous and strictly increasing. Further more assume that there exists an identity element $\tilde{e} \in \tilde{R}$ such that

$$x * \tilde{e} = x, \quad x \in A. \tag{1.3}$$

---

[*]The first author acknowledges research support from DST, Govt of India.

Further more every continuous, reducible and associative operation defined on an interval $A$ in the real line is communicative (Aczel (1966)). Characterization of distributions through binary operations is given in Muliere and Scarcini (1987) and Muliere and Prakasa Rao (2003). In Prakasa Rao (2004), the bivariate lack of memory property (Roy (2002)) is generalized and classes of bivariate probability distributions which include bivariate exponential, Weibull, Pareto distributions are characterized under binary associative operations.

Asha and Rejeesh (2007) characterized a subclass of the reversed generalized Pareto distribution (Castillo and Hadi (1995)) using the reversed lack of memory property. In this paper we generalize the reversed lack of memory property and characterize probability distributions using this property. In Section 2 the univariate reversed lack of memory property is generalized using the binary operation and a class of probability distributions which include subclass of reversed generalized Pareto distribution, power function and reflected Weibull is characterized. In Section 3 this property is extended to the bivariate set up the class of bivariate probability distributions which are characterized by this property is derived.

## 2. Univariate characterizations

Let $X$ be an absolutely continuous random variable with density, distribution function and survival function defined as $f(\cdot)$, $F(\cdot)$ and $R(\cdot)$ respectively. Let $a = \inf\{x : F(x) > 0\}$ and $b = \sup\{x : F(x) < 1\}$, then $(a, b)$, $-\infty \leq a < b \leq \infty$ is the interval of support of $X$. A random variable $X$ with support $(a, b)$ where $a < 0$ and $b \geq 0$ is said to have the reversed lack of memory property (RLMP) (Asha and Rejeesh (2007)), if

$$P(X \leq x | X \leq x + t) = P(X \leq 0 | X \leq t)$$

for all $x$ and $t$ such that $a < x \leq x + t \leq b$.
In terms of the distribution function we can write the RLMP as

$$F(x) \cdot F(t) = F(x + t) \cdot F(0) \tag{2.1}$$

where $a < x \leq x + t \leq b$ with $a < 0$.
Solving (2.1) (Aczel, 1966, p. 38) we have

$$F(x) = \begin{cases} e^{c(x-b)}, & x < b \\ 1, & x \geq b, \ c > 0. \end{cases} \tag{2.2}$$

This type 3 extreme value distribution also belongs the the reversed generalized Pareto distribution developed by Castillo and Hadi (1995) as a fatigue model that satisfy certain compatibility conditions arising out of physical and statistical conditions in fatigue studies.

We shall extend the RLMP using the binary operator $*$ as

$$P(X \leq x | X \leq x * t) = P(X \leq e | X \leq t) \tag{2.3}$$

for all $a < x < b, a < e, e \leq t \leq b, x * t \leq b$.

Here $e$ is specified as in (1.3). The general reducible continuous solution is (1.2).

In terms of the distribution function (2.3) can be written as

$$F(x) \cdot F(t) = F(x * t) \cdot F(e) \tag{2.4}$$

for all $a < x < b, a < e, e \leq t \leq b, x * t \leq b$.

We refer to (2.4) as generalized reversed lack of memory property (GRLMP).

Recently the reversed hazard function has become quite popular since it is very useful in the analysis of left censored data (see, Gupta and Hann (2001), Andersen et al. (1993)). The reversed hazard rate (RHR) function of $X$ is defined for $x > a$ as

$$\lambda(x) = \frac{d}{dx} \log F(x) = \frac{f(x)}{F(x)}.$$

It has been observed that there exists a relation between proportional reversed hazard class of distributions and the exponentiated class of distributions. The class of distributions $[F(\cdot)]^\alpha$, $\alpha > 0$ can be defined as the exponentiated class of distributions with base distributions $F(\cdot)$. In this case the reversed hazard rate function of $X$ satisfy $\lambda_Y(t) = \alpha \lambda_X(t)$ where $Y$ has a distribution $[F(\cdot)]^\alpha$. Observe that $F(x)$ in (2.2) is a proportional reversed hazards model.

In the next theorem we show that the continuous solution of (2.4) are generalized proportional reversed hazard (PRH) models.

**Theorem 2.1.** *The (continuous) solution of* (2.4) *is*

$$F(x) = \exp[c(g(x) - g(b))] \quad \text{with } c > 0 \text{ and } x \in (g^{-1}(a), b).$$

*where $g$ is a continuous and strictly monotone function.*

*Proof.* Substituting (1.2) in (2.4) we have

$$F(x) \cdot F(t) = F[g^{-1}(g(x) + g(t))] \cdot F(e). \tag{2.5}$$

Writing

$$s = g(x), \ u = g(t) \quad \text{and} \quad F \circ g^{-1} = H \tag{2.6}$$

(2.5) becomes

$$H(s) \cdot H(u) = H(s + u) \cdot H(g(e))$$

for all $g(a) < s < g(b)$ and $g(e) \leq u \leq g(b)$ which implies

$$G(s) \cdot G(u) = G(s + u) \quad \text{where} \quad G(s) = \frac{H(s)}{H(g(e))}$$

for all $g(a) < s < g(b)$ and $g(e) \leq u \leq g(b)$ which is the Cauchy functional equation and the solution is

$$G(s) = e^{cs}, \quad s \in (g(a), g(b)) \quad \text{and} \quad c > 0 \quad \text{(Aczel, 1966, p. 38)}.$$

From (2.6) we get

$$F(x) = e^{cg(x)} \cdot F(e), \quad x \in (g^{-1}(a), b). \tag{2.7}$$

Now taking $x = b$ we get $F(e) = e^{-cg(b)}$.
Thus (2.7) reduces to

$$F(x) = \begin{cases} e^{c[g(x)-g(b)]}, & x \in (g^{-1}(a), b) \\ 1, & x \geq b, \ c > 0. \end{cases}$$

Hence the theorem. □

If we particularize the operation $*$, we obtain different types of distributions.

**Example 1.** Subclass of reversed generalized Pareto distribution (Castillo and Hadi (1995)).

For $x * y = x + y$, we get $g(x) = x$, $x \in (-\infty, b)$ with $a = -\infty$, $b < \infty$ and $e = 0$.
The distribution function is now given by

$$F(x) = \begin{cases} \exp[c(x-b)], & x < b \\ 1, & x \geq b, \ c > 0 \end{cases}$$

**Example 2.** Power function distribution.

(i) For $x * y = xy$, we get $g(x) = \log x$, $x \in (0, b)$ with $e = 1$. In this case the distribution function is given by

$$F(x) = \begin{cases} (\frac{x}{b})^c, & 0 \leq x < b \\ 1, & x \geq b, \ c > 0. \end{cases}$$

which is the distribution function for the power function distribution.

(ii) If we take $x * y = x + y + xy$, we get $g(x) = \log(x+1)$, $x \in (-1, b)$, $b < \infty$ with $e = 0$. Here the distribution function is given by

$$F(x) = \begin{cases} (\frac{x+1}{b+1})^c, & -1 \leq x < b \\ 1, & x \geq b, \ c > 0. \end{cases}$$

which is the distribution function for the power function distribution in the support of $(-1, b), b < \infty$.

**Example 3.** Reflected Weibull distribution (Lai and Xie (2005)).

If we take $x * y = \sqrt{x^2 + y^2}$, we get $g(x) = -x^2$, $x \in (-\infty, 0)$, with $e = 0$. The corresponding distribution function is given by

$$F(x) = \begin{cases} e^{-cx^2}, & x < 0 \\ 1, & x \geq 0, \ c > 0 \end{cases}$$

which is the distribution function for a reflected Weibull distribution.

**Remark 2.1.** *If $X$ is a random variable in the support of $(g^{-1}(e), b)$ with $F(e) \neq 0$, then $X$ satisfies the GRLMP if and only if*

$$F_X(x) = \begin{cases} \exp[c[g(x) - g(b)]], & g^{-1}(e) \leq x \leq b \\ 1, & x \geq b, \ c > 0 \end{cases}$$

*which has a probability mass at $x = g^{-1}(e)$.*

Now we try to evolve the above concept to the higher dimensions. One of the main problems associated with such an attempt is that there is no unique way of evolution. We here consider the definition of bivariate reversed lack of memory property (Asha and Rejeesh (2007)) to extend the concept of univariate GRLMP to the bivariate case and derive bivariate models characterized by the respective property. Since multivariate derivations are rather straight forward extensions it is excluded.

## 3. Bivariate Extensions

Consider a random vector $X = (X_1, X_2)$ in the two-dimensional space with joint distribution function $F(x_1, x_2) = P(X_1 \leq x_1, X_2 \leq x_2)$ in the support of $(a, b)^2$ where

$$a = \inf\{x_i | F_i(x_i) > 0\} \quad \text{and} \quad b = \sup\{x_i | F_i(x_i) < 1\},$$

with $F(x_1, b)$ and $F(b, x_2)$ as the marginals of $X_i$, $i = 1, 2$.
Then a direct extension of RLMP is

$$F(x_1 + t_1, x_2 + t_2) \cdot F(0, 0) = F(x_1, x_2) \cdot F(t_1, t_2) \tag{3.1}$$

for all $x_i$ and $t_i$ such that $a < x_i \leq x_i + t_i \leq b$, $a < 0$, $i = 1, 2$.
The only solution for (3.1) is (Aczel (1966))

$$F(x_1, x_2) = e^{\lambda_1(x_1 - b) + \lambda_2(x_2 - b)}, a < x_i < b, \ \lambda_i > 0, \ i = 1, 2.$$

We consider analogous equation of (3.1) for the associative binary operator $*$ given by

$$F(x_1 * t_1, x_2 * t_2) \cdot F(e, e) = F(x_1, x_2) \cdot F(t_1, t_2) \tag{3.2}$$

for all $a < x_i < b$, $e \leq t_i \leq b$, $a < e$, $x_i * t_i \leq b$, $i = 1, 2$.
Combining (1.2) and (3.2) we have

$$F[g^{-1}(g(x_1) + g(t_1)), g^{-1}(g(x_2) + g(t_2))] \cdot F(e, e) = F(x_1, x_2) \cdot F(t_1, t_2). \tag{3.3}$$

Writing $s_i = g(x_i)$, $u_i = g(t_i)$, $i = 1, 2$ and $F(g^{-1}(\cdot), g^{-1}(\cdot)) = H(\cdot, \cdot)$, (3.3) becomes

$$H(s_1 + u_1, s_2 + u_2) \cdot H(g(e), g(e)) = H(s_1, s_2) \cdot H(u_1, u_2)$$

for all $g(a) < s_i < g(b)$ and $g(e) \leq u_i \leq g(b)$, $i = 1, 2$.
Thus

$$G(s_1 + u_1, s_2 + u_2) = G(s_1, s_2) \cdot G(u_1, u_2) \text{ where } G(s_1, s_2) = \frac{H(s_1, s_2)}{H(g(e), g(e))}$$

for all $g(a) < s_i < g(b)$ and $g(e) \leq u_i \leq g(b)$, $i = 1, 2$.
The solution to the above Cauchy functional equation is (Aczel (1966))

$$G(s_1, s_2) = e^{\lambda_1 s_1 + \lambda_2 s_2}, \lambda_1, \lambda_2 > 0$$

which implies

$$F(x_1, x_2) = e^{\lambda_1 [g(x_1) - g(b)] + \lambda_2 [g(x_2) - g(b)]}$$

for all $x_i \in (g^{-1}(a), b)$, $\lambda_i > 0$, $i = 1, 2$.
Thus, $X_1$ and $X_2$ are independent with marginal distribution functions specified by

$$F_{X_1}(x_1) = e^{\lambda_1 [g(x_1) - g(b)]} \quad \text{and} \quad F_{X_2}(x_2) = e^{\lambda_2 [g(x_2) - g(b)]}$$

We now consider a meaningful way of extending the RLMP specified by (2.1) to the bivariate case as

$$F(x_1 + t, x_2 + t) \cdot F(0, 0) = F(x_1, x_2) \cdot F(t, t) \tag{3.4}$$

for all $a < x_i \leq x_i + t \leq b$, $i = 1, 2$ with $a < 0$.
Then the unique solution of (3.4) among distribution function is (Asha and Rejeesh (2007))

$$F(x_1, x_2) = e^{\lambda_1 (x_1 - b) + \lambda_2 (x_2 - b) + \lambda_{12} \max(x_1 - b, x_2 - b)}$$

for all $a < x_i < b$; $\lambda_i, \lambda_{12} > 0$; $i = 1, 2$.
Generalizing equation (3.4) using the operation $*$. We define the generalized bivariate reversed lack of memory property (GBRLMP) as

$$F(x_1 * t, x_2 * t) \cdot F(e, e) = F(x_1, x_2) \cdot F(t, t) \tag{3.5}$$

for all $a < x_i < b$, $e \leq t \leq b$, $x_i * t \leq b$, $a < e$, $i = 1, 2$.

**Theorem 3.1.** *Let* $F(x_1 * t, x_2 * t) \cdot F(e, e) = F(x_1, x_2) \cdot F(t, t)$ *for all* $a < x_i < b$, $e \leq t \leq b$, $x_i * t_i \leq b$, $a < e$; $i = 1, 2$ *and*

$$F_i(x * t) \cdot F_i(e) = F_i(x) \cdot F_i(t) \tag{3.6}$$

*for all* $a < x < b$, $e \leq t \leq b$, $x * t \leq b$, $a < e$; $i = 1, 2$. *Then the continuous solution of the equations* (3.5) *and* (3.6) *is*

$$F(x_1, x_2) = \exp\{\lambda_1 [g(x_1) - g(b)] + \lambda_2 [g(x_2) - g(b)] \\ + \lambda_{12} \max[g(x_1) - g(b), g(x_2) - g(b)],$$

$g^{-1}(a) < x_i < b$; $\lambda_i > 0$; $i = 1, 2$, $\lambda_{12} \geq 0$.

*Proof.* In (3.5) let $x_1 = x_2 = x$, then

$$F(x * t, x * t) \cdot F(e, e) = F(x, x) \cdot F(t, t). \tag{3.7}$$

Combining (1.2) and (3.7) we get

$$F[g^{-1}(g(x) + g(t)), g^{-1}(g(x) + g(t))] \cdot F(e, e) = F(x, x) \cdot F(t, t).$$

Taking $g(x) = s$, $g(t) = u$ and $F(g^{-1}(\cdot), g^{-1}(\cdot)) = H(\cdot, \cdot)$. Thus (3.7) reduces to the Cauchy functional equation

$$G(s+u) = G(s) \cdot G(u) \qquad (3.8)$$

for all $g(a) < s < g(b)$ and $g(e) \leq u \leq g(b)$ where $G(s) = \frac{H(s,s)}{H(g(e),g(e))}$.
Solving for (3.8) (Aczel, 1966, p. 38) we have

$$G(s) = e^{ks} \text{ for } k > 0.$$

Thus we get $\frac{F(x,x)}{F(e,e)} = e^{kg(x)}$, $k > 0$ which implies

$$F(x,x) = e^{k[g(x)-g(b)]}, \ x \in (g^{-1}(a), b), k > 0. \qquad (3.9)$$

Hence

$$F(x_1 * t, b * t) = \frac{F(x_1, b) \cdot F(t, t)}{F(e, e)}$$
$$= \frac{F_1(x_1) \cdot F(t, t)}{F(e, e)}$$

From (3.9) we get

$$F(x_1 * t, b * t) = \frac{e^{\lambda_1[g(x_1)-g(b)]+k[g(t)-g(b)]}}{e^{-kg(b)}}$$
$$= e^{\lambda_1[g(x_1)-g(b)]+kg(t)}. \qquad (3.10)$$

Now let $x_1 * t = s$ and $b * t = u$, then by (1.2) we get $g(x_1) = g(s) - g(t)$ and $g(b) = g(u) - g(t)$. Thus (3.10) becomes

$$F(s,u) = e^{\lambda_1[g(s)-g(u)]+k[g(u)-g(b)]}; \quad u \geq s.$$

Arguing similarly we can prove that

$$F(s,u) = e^{\lambda_2[g(u)-g(s)]+k[g(s)-g(b)]}; \quad u \leq s.$$

Writing $\lambda_{12} = k - \lambda_1 - \lambda_2$ and rearranging we obtain,

$$F(s,u) = \exp\{\lambda_1[g(s) - g(b)] + \lambda_2[g(u) - g(b)]$$
$$+ \lambda_{12} \max[g(s) - g(b), g(u) - g(b)]$$

or

$$F(x_1, x_2) = e^{\lambda_1[g(x_1)-g(b)]+\lambda_2[g(x_2)-g(b)]+\lambda_{12}\max[g(x_1)-g(b),g(x_2)-g(b)]},$$
$$g^{-1}(a) < x_i < b; \lambda_i > 0; i = 1,2 \text{ and } \lambda_{12} \geq 0. \qquad (3.11)$$

□

**Remark 3.1.** *It can be seen that a necessary and sufficient condition for $X_1$ and $X_2$ to be independent is $\lambda_{12} = 0$.*

**Remark 3.2.** *If $(X_1, X_2)$ is distributed as in Theorem 3.1, then the distribution function of $\max(X_1, X_2)$ has the following form:*

$$\begin{aligned}P[\max(X_1, X_2) < s] &= P[X_1 < s, X_2 < s] \\ &= F(s, s) \\ &= e^{(\lambda_1 + \lambda_2 + \lambda_{12})(g(s) - g(b))}.\end{aligned}$$

*Hence, the distribution of $\max(X_1, X_2)$ has the same form (with different parameter) as the marginal distribution of $X_1$ and $X_2$.*

**Theorem 3.2.** *Let $(X_1, X_2)$ be a bivariate random variable with joint distribution function $F(x_1, x_2)$. Then F is distributed as in (3.11) if there exist an independent and identically distributed random variables $U, V, W$ whose marginal distributions are given by $F(x) = e^{\lambda[g(x) - g(b)]}$ such that $X_1 = \max(U, W)$, $X_2 = \max(V, W)$.*

*Proof.*

$$\begin{aligned}F(x_1, x_2) &= P(X_1 < x_1, X_2 < x_2) \\ &= P(U < x_1, W < x_1, V < x_2, W < x_2) \\ &= P(U < x_1) \cdot P(V < x_2) \cdot P(W < \max(x_1, x_2)) \\ &= e^{\lambda_1[g(x_1) - g(b)] + \lambda_2[g(x_2) - g(b)] + \lambda_{12}\max[g(x_1) - g(b), g(x_2) - g(b)]}.\end{aligned}$$

$\square$

The distribution function $F(x_1, x_2)$ specified in Theorem 3.1 is not an absolutely continuous distribution in the support of $(g^{-1}(e), b)^2$ with $F(e, e) \neq 0$. In this case we can write

$$F(x_1, x_2) = F_d + F_c + F_{ac}$$

where

$$F_d(x_1, x_2) = \begin{cases} e^{\lambda_i(e - g(b))}; & x_i = g^{-1}(e),\ x_{3-i} = b,\ i = 1, 2. \\ e^{k(e - g(b))}; & x_1 = x_2 = g^{-1}(e),\ k = \lambda_1 + \lambda_2 + \lambda_{12}. \end{cases}$$

$$F_c(x_1, x_2) = e^{k[g(x) - g(b)]},\ x_1 = x_2 = x,\ g^{-1}(e) < x < b$$

and

$$F_{ac}(x_1, x_2) = e^{\lambda_1[g(x_1) - g(b)] + \lambda_2[g(x_2) - g(b)] + \lambda_{12}\max[g(x_1) - g(b), g(x_2) - g(b)]}$$

$g^{-1}(e) < x_i < b$; $\lambda_i > 0$; $i = 1, 2$; $\lambda_{12} \geq 0$. Few members of (3.11) are listed in table 1

Table 1.

| No | $x*y$ | $e$ | $g(x)$ | Distribution |
|----|-------|-----|--------|--------------|
| 1 | $x+y$ | 0 | $g(x)=x$, $x\in(-\infty,b), b<\infty$ | Bivariate Type 3 extreme value distribution $F(x_1,x_2)=\exp[\lambda_1(x_1-b)+\lambda_2(x_2-b)+\lambda_{12}\max(x_1-b,x_2-b)]$ $-\infty<x_i<b; \lambda_i>0, i=1,2; \lambda_{12}\geq 0$. |
| 2 | $xy$ | 1 | $g(x)=\log x$, $x\in(0,b), b<\infty$ | Bivariate power function distribution $F(x_1,x_2)=(\frac{x_1}{b})^{\lambda_1}(\frac{x_2}{b})^{\lambda_2}\cdot e^{\lambda_{12}\max[\log(\frac{x_1}{b}),(\frac{x_2}{b})]}$ $0<x_i<b; \lambda_i>0, i=1,2; \lambda_{12}\geq 0$. |
| 3 | $x+y+xy$ | 0 | $g(x)=\log(x+1)$, $x\in(-1,b), b<\infty$ | Bivariate power function distribution $F(x_1,x_2)=(\frac{x_1+1}{b+1})^{\lambda_1}(\frac{x_2+1}{b+1})^{\lambda_2}e^{\lambda_{12}\max[\log(\frac{x_1+1}{b+1}),\log(\frac{x_2+1}{b+1})]}$ $-1<x_i<b; \lambda_i>0, i=1,2; \lambda_{12}\geq 0$. |
| 4 | $\sqrt{x^2+y^2}$ | 0 | $g(x)=-x^2$, $x\in(-\infty,0)$ | Bivariate reflected Weibull distribution $F(x_1,x_2)=\exp[-\lambda_1 x_1^2-\lambda_2 x_2^2+\lambda_{12}\max(-x_1^2,-x_2^2)]$ $-\infty<x_i<0; \lambda_i>0, i=1,2; \lambda_{12}\geq 0$. |